\documentclass[12pt,reqno,twoside]{amsart}
\usepackage{amsmath}
\usepackage{amssymb}
\usepackage{epsfig}
\usepackage{color}

\numberwithin{equation}{section}

\title[Friendly measures and singular vectors]{Friendly measures,
homogeneous flows and singular vectors}

\author{Dmitry Kleinbock}

\address{Brandeis University, Waltham MA
02454-9110 {\tt kleinboc@brandeis.edu}}

\author{Barak Weiss}

\address{Ben Gurion University, Be'er Sheva, Israel 84105
{\tt barakw@math.bgu.ac.il}}

\newif\ifdraft\drafttrue

\draftfalse

\long\def\comdima#1{\ifdraft{\marginpar{$\bullet$ \it
#1}}\else\ignorespaces\fi}
\long\def\combarak#1{\ifdraft{{\sb #1} }\else\ignorespaces\fi}

\newcommand{\Q}{{\mathbb {Q}}}

\newcommand{\compose}{{\circ}}

\newcommand{\R}{{\mathbb{R}}}

\newcommand{\Z}{{\mathbb{Z}}}

\newcommand{\E}{{\mathbf{e}}}

\newcommand{\N}{{\mathbb{N}}}

\newcommand{\cl}{\overline}

\newcommand{\ww}{{\bf{w}}}

\newcommand{\vf}{{\bf{f}}}

\newcommand{\vv}{{\bf{v}}}

\newcommand{\vr}{{\bf{r}}}

\newcommand{\ve}{{\bf{e}}}

\newcommand{\vw}{{\bf{w}}}

\newcommand{\GL}{\operatorname{GL}}

\newcommand{\SL}{\operatorname{SL}}

\newcommand{\ggm}{G/\Gamma}

\newcommand{\diag}{{\rm diag}}

\newcommand{\df}{{\, \stackrel{\mathrm{def}}{=}\, }}

\newcommand{\x}{{\bf x}}

\newcommand{\z}{{\bf z}}

\newcommand{\p}{{\bf p}}

\newcommand{\vn}{{\bf n}}

\newcommand{\til}{\widetilde}

\newcommand{\supp}{{\rm supp}}

\newcommand{\ff}{{\mathbf{f}}}

\newcommand{\pp}{{\mathbf{p}}}

\newcommand{\sm}{\smallsetminus}

\newcommand{\vre}{\varepsilon}

\font\sb = cmbx8 scaled \magstep0

\newcommand\ba{badly approximable}

\newcommand\di{diophantine}

\newcommand{\Sing}{{\mathrm{Sing}}}

\newcommand {\equ}[1]     {\eqref{#1}}

\newcommand {\ignore}[1]  {}

\newtheorem{thm}{Theorem}[section]

\newtheorem{lem}[thm]{Lemma}

\newtheorem{prop}[thm]{Proposition}

\newtheorem{cor}[thm]{Corollary}

\newtheorem{remark}[thm]{Remark}

\begin{document}

\ignore{

\begin{abstract}

\end{abstract}

}


\maketitle


\section{Introduction}

The theory of diophantine approximation studies how well $\x  = (x_1,
\ldots, x_n) \in \R^n$
can be approximated by $(p_1/q_1, \ldots, p_n/q_n)
\in \Q^n$ of a given `complexity', where this complexity is
usually measured by the quantity $\mathrm{lcm}(q_1, \ldots, q_{n})$.
Thus one is
interested in minimizing the difference, in a suitable sense, between
$q\x$ and a vector  $\pp$, where $\pp \in \Z^n$ and $q \in
\N$, with a given upper bound on $q$. Often one finds that certain
approximation problems admit a solution for almost every $\x$, while
others admit a solution for almost no $\x$; one is then interested in
understanding whether the typical properties remain typical when additional
restrictions are placed on $\x$.

As an example, consider the notion of a {\em singular\/} vector,
introduced by  A.\ Khintchine in the 1920s (see \cite{Khintchine,
Cassels}). 
Say that $\x$ is
{\em singular\/} if for any $\delta>0$ there is $T_0$ such that for all
$T \geq T_0$ one
can find $\pp \in \Z^n$ and $q \in \N$ with
\begin{equation}
\label{eq: singular}
\|q\x - \pp\| < \frac{\delta}{T^{1/n}}  \ \ \ \mathrm{and}
\ \ \ q  < T\,.
\end{equation}
a dual form,
Clearly this definition is independent on the choice of the norm.
Note also that by Dirichlet's Theorem, when $\delta > 1$ and
$\|\cdot\|$ is chosen to
be the supremum norm, the system
\equ{eq: singular} has
a nonzero integer solution for any $T > 1$. Thus singular vectors are
often referred to as
those for which  Dirichlet's theorem can be infinitely improved.

Let us say that $\x$ is  {\em totally irrational} if
$1, x_1, \ldots, x_n$ are linearly independent over $\Q$. It is not
hard to see
that vectors which are
not  totally irrational  are singular, and that the converse is true
    for
$n=1$. \comdima{I guess you are right, one needs transference.} However,
for  $n>1$  Khintchine \cite{Khintchine} proved the
existence of  totally irrational singular vectors.
On the other hand it is straightforward  to verify
\cite[Ch.\,V, \S7]{Cassels} that
Lebesgue measure of the set of singular vectors is zero.
In the late 1960s H.\ Davenport and W.\ Schmidt
showed \cite{Davenport-Schmidt} that
$\x \in \R^2$ of the form $\x = (t, t^2)$ is not
singular for Lebesgue-a.e.\ $t\in\R$. This was later extended to
certain classes of smooth curves and higher-dimensional submanifolds
of $\R^n$ by R.~Baker
\cite{Baker JLMS, Baker MPCPS} and
M.\ Dodson, B.\ Rynne and  J.\ Vickers \cite{DRV1} respectively; see
\S\ref{inher}
for precise statement of their results.

In this paper
we consider several generalizations of the notion of a singular vector.
Namely, following \cite{dima duke, dima tata, PV-bad}, we attach
different weights
to different components of $\x$ by means of the
{\em $\vr$-quasinorm\/}
\begin{equation*}
\label{eq: dfn quasinorm}
\| \x \|_{\vr} \df \max_{i=1, \ldots, d} |x_i|^{1/r_i}\,,
\end{equation*}
where
\begin{equation}
\label{eq: setting r}
\vr =
(r_1,\dots,r_{n}) \quad \mathrm{with} \quad
r_i > 0\quad\text{and}\quad\sum_{i=1}^n r_i = 1\,.
\end{equation}
Then  say that
    $\x$ is
{\em $\vr$-singular\/} if for any $\delta>0$ there is $T_0$
such that for all
$T \ge T_0$  one
can find $\pp \in \Z^n$ and $q \in \N$ with \comdima{I kept your unorthodox
form of the inequality, maybe we indeed need to comment
on its difference from the conventional one.}
\begin{equation}
\label{eq: r-T-singular}
\|q\x - \pp\|_\vr < \frac{\delta}{T}  \ \ \ \mathrm{and}
\ \ q  < \delta T\,.
\end{equation}
Further, for $\vr$ as above and an unbounded subset $\mathcal{T}$ of $\R_{\ge
1}$,  say that
    $\x$ is
{\em $\vr$-singular along\/} $\mathcal{T}$ if for any $\delta>0$ there is $T_0$
such that for all
$T \in \mathcal{T} \cap [T_0,\infty)$  one
can find $\pp \in \Z^n$ and $q \in \N$ satisfying \equ{eq: r-T-singular}.
We will denote the set of $\vr$-singular (along $\mathcal{T}$) vectors
by $\Sing(\vr)$ and
$\Sing\left(\vr,\mathcal{T}\right)$ respectively.

It is clear that $\x$ is
singular if and only if
$\x\in\Sing(\vn)$, where
    $$\vn \df (1/n,\dots,1/n)\,,$$ the vector assigning equal weights to
each coordinate, and that  $\Sing(\vr) = \Sing\left(\vr,\R_{\ge 1}\right)$
is contained  in $\Sing\left(\vr,\mathcal{T}\right)$
for any
    $\mathcal{T}\subset\R_{\ge 1}$. However an elementary modification
of the proof
in \cite[Ch.\,V, \S7]{Cassels}  shows \comdima{This is completely
straightforward,
so there is no point to use ergodicity/mixing; let us reserve it for the
   DI paper.}
that Lebesgue measure of $\Sing\left(\vr,\mathcal{T}\right)$ is zero for any
$\vr$ as in \equ{eq: setting r} and any unbounded $\mathcal{T}$. In this paper
we consider the class of {\em friendly\/} measures on $\R^n$,
originally introduced in
\cite{friendly} and described in detail in \S\ref{inher}, and
prove

\begin{thm}\label{cor: friendly}
If $\mu$ is a friendly measure on $\R^n$, then for any $\vr$ as in
\equ{eq: setting r} and any unbounded $\mathcal{T}$,
$\mu\big(\Sing(\vr,\mathcal{T}) \big)=0$.
\end{thm}

A special case of this theorem, with $\vr = \vn$ and $\mathcal{T} =
\R_{\ge 1}$,
was announced in \cite{friendly}.

\medskip

    The class of
friendly measures includes
Hausdorff measures supported on various self-similar sets such as
Cantor's ternary set, Koch snowflake, Sierpinski gasket, etc.
It also includes  volume measures on smooth nondegenerate
submanifolds of $\R^n$.
We recall that $M\subset \R^n$ is called {\em nondegenerate\/} if it
is parameterized by
a smooth map $\vf$ from an open subset $U$ of $\R^d$ to $\R^n$ such
that for Lebesgue-a.e.\
$\x\in U$ there exists $\ell\in\N$ such that partial derivatives of
$\vf$ at $\x$
up to order $\ell$ span $\R^n$. If $\vf$ is real analytic and $U$ is connected,
   the latter
condition is equivalent to $\vf(U)$ not
being contained in a proper affine hyperplane of $\R^n$.
Thus Theorem \ref{cor: friendly} significantly generalizes the
aforementioned result of \cite{Davenport-Schmidt} about the curve
$\{(t,t^2): t \in \R \}$, as
well as additional results obtained by several authors.

Note that it is not hard to construct a friendly measure whose support does not
contain any singular vectors at all; for example, the set of \ba\
vectors supports
friendly measures of arbitrarily small codimension \cite{bad, Urbanski
- d=1}. The
situation is however different for volume measures on real analytic
nondegenerate manifolds.
Namely, we prove

\begin{thm}
\label{thm: main, dim2}
Let ${M} \subset \R^n$ be a real analytic
submanifold of
dimension at least $2$  which is not contained in any proper rational affine
hyperplane of $\R^n$,
and let $\vr$ be as in \equ{eq: setting r}.
Then there exists a totally irrational $\x \in
{M}\cap {\Sing}(\vr)$.
\end{thm}

Our approach to Theorem \ref{cor: friendly} is modelled
on \cite{KM, friendly}: in \S\ref{dynam} we translate
the aforementioned \di\
properties of $\x\in \R^n $ into dynamical properties of
certain trajectories in the homogeneous space $\ggm$,
where
\begin{equation}
\label{eq: defn ggm}
G=\SL(n+1, \R)\quad \text{and} \quad\Gamma = \SL(n+1, \Z)\,.
\end{equation}
Namely, we show
    (Proposition
\ref{prop: dynamical interpretation}) that $\x\in \R^n $ is
$\vr$-singular along  $\mathcal{T}$ if and only if
the corresponding trajectory  leaves every compact subset of $\ggm$.
To control  the measure of points with divergent trajectories we
employ quantitative nondivergence
estimates from \cite{friendly}, described in detail in \S\ref{quant}.
Theorem \ref{thm: main, dim2} is proved in \S\ref{cassels}; the argument
    is a modification of the proof of \cite[Thm.\ 5.2]{cassels},
and is based on ideas going back to Khintchine \cite{Khintchine}.
\ignore{The last section of the paper contains
some open questions and indications of more general results provable
by our methods.
In particular, much can be said about vectors $\x$ for which there
exists $\delta < 1$
and $T_0 > 0$ such that the system \equ{eq: singular},
or, more generally \equ{eq: r-T-singular}, has an integer solution
for any $T > T_0$
(i.e.\ those for which $1$ in Dirichlet's Theorem
can be improved  to $\delta$). More detail will follow in a forthcoming sequel
\cite{KW} to this paper.}

\medskip

{\bf Acknowledgements:}
This research was supported by BSF grant
2000247 and NSF grant DMS-0239463.
We are grateful to
the Max Planck Institute for its hospitality during July 2004, and in
particular to Sergiy Kolyada, who organized the activity at MPI. We are
also grateful to Roger
Baker for useful discussions.



\section{Dynamical interpretation of singular vectors}
\label{dynam}

Let $G$ and $\Gamma$ be as in \equ{eq: defn ggm}, and denote by $\pi$
the quotient map from $G$ onto $\ggm$.
    $G$ acts on $\ggm$ by left translations via the rule $g \pi(h) =
\pi(gh)$, $g,h\in G$.
Define
\begin{equation*}
\label{eq: defn tau}
\tau(\x) \stackrel{\mathrm{def}}{=} \left(
\begin{array}{ccccc}
I_{n} & \x \\ 0 & 1
\end{array}
\right), \ \ \ \bar{\tau} \stackrel{\mathrm{def}}{=} \pi \compose \tau\,,
\end{equation*}
where $I_{n}$ stands for the $n\times n$ identity matrix.
Then, given $\vr$
as in \equ{eq: setting r}, consider the  one-parameter subgroup
$\{g_t^{(\vr)}\}$ of $G$
given by
\begin{equation*}
\label{eq: new defn g_t}
g_t^{(\vr)} \df \diag(e^{r_1t}, \ldots,
e^{r_{n} t}, e^{-t})\,.
\end{equation*}

Recall that $\ggm$ is noncompact.
For an unbounded subset $A$ of $\R_+$ and $x\in\ggm$, say that a trajectory
$\{g^{(\vr)}_t x: t \in A\}$ is {\em divergent\/}
if the map $A\to \ggm$, $t \mapsto g^{(\vr)}_t x$, is proper; that is,
    for any compact $K \subset \ggm$ there exists
$t_0$ such that  $g^{(\vr)}_t \notin K$ for all $t \in A \cap [t_0,\infty)$.

It was proved in \cite[Proposition 2.12]{Dani} that $\x$ is singular
if and only if
the trajectory $\{g^{(\vn)}_t \bar{\tau}(\x)
: t \geq 0\}$  in $\ggm$ is divergent, and
in \cite[Theorem 7.4]{dima duke} that $\x$ is $\vr$-singular if and only if
the trajectory $\{g^{(\vr)}_t \bar{\tau}(\x)
: t \geq 0\}$ is divergent. We generalize this correspondence one step further:

\begin{prop}
\label{prop: dynamical interpretation}
    $\x$ is
$\vr$-singular along  $\,\mathcal{T}$ if and only if the trajectory
\begin{equation}
\label{traj}
\{g^{(\vr)}_t \bar{\tau}(\x)
: t \in\log \mathcal{T}\}\subset \ggm
\end{equation}
   is divergent, where $\log \mathcal{T} \df \{\log T : T\in
    \mathcal{T}\}$.
\end{prop}

To prove Proposition \ref{prop: dynamical interpretation}, we need
an explicit description of compact subsets of $\ggm$.
Since $\Gamma$ is the stabilizer of $\Z^{n+1}$ under the action of $G$ on
the set of lattices in $\R^{n+1}$, $\ggm$ can be identified with $G\Z^{n+1}$,
that is, with the set of all unimodular lattices in $\R^{n+1}$.
Fix a norm $\|\cdot\|$ on  $\R^{n+1}$, and for $\varepsilon>0$ let
\begin{equation}
\label{k epsilon}
\begin{split}
K_{{\varepsilon}} &\stackrel{\mathrm{def}}{=} \pi\big(\big\{ g \in G : \Vert g
\vv \Vert \geq {\varepsilon} \quad \forall\, \vv \in
\Z^{n+1} \sm \{ 0 \}\big\}\big);
\end{split}
\end{equation}
i.e., $K_{ \varepsilon}$ is the collection of all unimodular
lattices in $\R^{n+1}$ which contain no nonzero vector with norm less than
$\varepsilon$.
By Mahler's compactness criterion (see e.g.\  \cite[Chapter
10]{Rag}), each $K_{{\varepsilon}}$ is compact,
and for each compact $K \subset \ggm$  there is $\vre>0$
such that $K \subset K_{ \vre}$.

Now take $\vr$ as in \equ{eq: setting r} and write
\comdima{Changed notation, hope you won't object}
$$\bar r = \min_{1 \leq i \leq n} r_i\,.%
$$



\begin{lem}\label{lem: correspondence}
Let $\|\cdot\|$ be the supremum norm, let $\vre$ and $t$ be positive
numbers with
$e^{\bar r  t} \geq \vre,$ and denote
$T = e^t$. Then
\equ{eq: r-T-singular} with $\delta = \vre^{1/\bar r }$
implies $g^{(\vr)}_t \bar{\tau}(\x) \notin
K_{\vre}$, which in turn implies \equ{eq: r-T-singular} with $\delta = \vre.$
\end{lem}

\begin{proof}
Suppose \equ{eq: r-T-singular}
holds with $\delta =
\vre^{1/\bar r }$ and with $\pp \in \Z^n, \, q
\in \N$. This implies that
$$e^{-t}q = q/T < \delta < \vre$$
and for $i=1, \ldots, d$,
$$e^{r_it} |
qx_i- p_i| < e^{r_it} \delta^{r_i}/T^{r_i} = \vre^{r_i/\bar r } \leq \vre\,.$$
   From this one concludes that for  $\vv = (-\pp, q) \in \Z^{n+1} \sm \{0\}$,
$$\|g_t^{(\vr)} {\tau}(\x) \vv \| = \max \left\{e^{-t}q,\, e^{r_1t} |
-p_1 +  qx_1|, \ldots,  e^{r_{n}t} |
-p_{n} +  qx_{n}|\right  \} < \vre\,,$$
so $g_t^{(\vr)} \bar{\tau}(\x) \notin K_{\vre}.$ The proof of the
second implication is
similar and is omitted.
\end{proof}

\begin{proof}[Proof of  Proposition~\ref{prop: dynamical interpretation}]
By the preceding lemma, $\x$ is
$\vr$-singular along  $\mathcal{T}$ if and only if  for any
$\varepsilon>0$ there is $T_0$
such that  $g^{(\vr)}_t \bar{\tau}(\x) \notin
K_{\vre}$ whenever $e^t \in \mathcal{T} \cap [T_0,\infty)$. The latter,
in view of  Mahler's compactness criterion,
    is equivalent to the fact that the trajectory \equ{traj}
   eventually leaves every compact subset of $\ggm$.
\end{proof}

As an application of this dynamical approach, we can state a condition
on a measure
$\mu$ on $\R^n$ guaranteeing that it assigns measure zero to singular vectors.

\begin{prop}
\label{prop: measure criterion} Let
$\mu$ be a measure on $\R^n$, and
suppose that  for $\mu$-a.e.\ $\x_0\in\R^n$ there is a ball  $B$
centered at $\x_0$
such that $\forall\,\delta > 0$  there exist $ \varepsilon > 0$ and a
sequence $t_k\to\infty,\
t_k\in\log\mathcal{T}$, with
\begin{equation}
\label{eq: ball}
\mu\big(\{\x\in B : g^{(\vr)}_{t_k}\bar{\tau}(\x)
\notin K_{{\varepsilon}}\}\big)  < \delta\text{ for every }k\,.\ \ \
\end{equation}
Then $\mu\big(\Sing(\vr,\mathcal{T}) \big)=0$.
\end{prop}

\begin{proof}
Indeed, if we take $\{t_k\}$ as above and let
$$
B_\varepsilon \df \bigcup_{N = 1}^\infty \bigcap_{k =
N}^\infty\{\x\in B : g^{(\vr)}_{t_k}\bar{\tau}(\x)
\notin K_{{\varepsilon}}\}\,,
$$ then \equ{eq: ball} implies that $\mu(B_\varepsilon) \le \delta$.
But the set
$\bigcap_\varepsilon B_\varepsilon$, in view of  Mahler's compactness
criterion,
    coincides with $$\{\x\in B : g^{(\vr)}_{t_k}\bar{\tau}(\x)
\text{ is divergent}\}\,,$$ and therefore  has measure zero.
\end{proof}

\begin{remark} \rm \comdima{The rest of the discussion is postponed
until the DI paper; instead I added
a remark on sup vs. Euclidean norms.}
Note that even though the definition of $K_{{\varepsilon}}$
depends on the choice of the norm $\|\cdot\|$ in \equ{k epsilon}, the
assumption
of Proposition \ref{prop: measure criterion} is clearly independent
of this choice.
   Thus without loss
of generality we may, and will, fix a Euclidean structure on $\R^{n+1}$
and choose  $\|\cdot\|$ to be the Euclidean norm.
\end{remark}

\section{The inheritance problem and friendly measures}
\label{inher}
In general, given a property which holds for a typical point in $\R^n$, it is
natural to inquire for which subsets (e.g.\ submanifolds, self-similar
`fractals', etc.) a typical point on the subset also
satisfies the property. The prototype for such an `inheritance
question' was the
famous conjecture of K.~Mahler from the 1930s,
settled three decades later by V.~Sprind\v zuk, which led to the
theory of diophantine
approximation on manifolds.

As was mentioned in the introduction, the first inheritance result related
    to the notion of singular vectors is due to Davenport and Schmidt:
they showed  that almost no points (with respect to the smooth measure class)
on $M$ are singular, where
\begin{itemize}
\item[(a)]  \quad  $M = \{(t, t^2): t \in \R\} \subset \R^2$
\cite[Thm.\ 3]{Davenport-Schmidt}.

\end{itemize}

Later, R.\ Baker and Y.\ Bugeaud proved that almost no points
on  $M$ are singular when:

\begin{itemize}

\item[(b)]  \quad $M = \{(t, t^2, t^3): t \in \R\} \subset \R^3$
\cite[Thm.\ 2]{Baker JLMS};

\item[(c)] \quad
$M= \{(t, \dots, t^n): t \in \R\} \subset \R^n$ \cite[Thm.\
7]{Bugeaud};

\item[(d)]  \quad $M \subset \R^2$ is a curve with continuous third
derivatives and
non-vanishing curvature almost everywhere \cite[Thm.\ 2]{Baker MPCPS}.

\end{itemize}


The only other  paper on this topic known to us is \cite{DRV2}, where
it was proved
that  almost all points on a  $C^3$ submanifold $M$ of $\R^n$  are
not singular if
\begin{itemize}

\item[(e)] \quad $M$ has `two-dimensional definite curvature almost
everywhere'.

\end{itemize}

The latter condition
requires the dimension of  $M$  to be at least $2$ (see \cite{DRV1}
for more detail).

The curve in (a--c) was first
studied by Mahler in connection with questions
about approximation of real numbers by algebraic numbers.
Note that all of the above examples are special case of nondegenerate
manifolds defined in the introduction
(see \cite[Remark 6.3]{KM} for a discussion of the relation between
nondegeneracy of $M$
and the
   conditions of \cite{DRV1, DRV2}.)
\ignore{
  However, it had not been
previously established
    that almost all points on $M$ are
    not singular when
$M$ is a nondegenerate submanifold of $\R^n$,
even for the special case
\begin{equation}
\label{mahler}
\end{equation}
which is the subject of the original
Mahler's Conjecture.
}

\medskip

A more general framework for discussing the inheritance problem is to
recast it in terms of measures. That is, given a property
$\mathcal{P}$ which holds
for Lebesgue-a.e. $\x \in \R^n$, one wants to describe measures $\mu$
such that $\mathcal{P}$ also holds for $\mu$-a.e.\ $\x$. Let us recall certain
properties of a measure on $\R^n$, introduced in
\cite{friendly}.

Suppose $\mu$ is a locally finite 
Borel measure on $\R^n$. Let
$B(\x,r)$ denote the
open ball of radius $r$ centered at
$\x$. Suppose $U\subset \R^n$ is open. We say that $\mu$ is
{\em $D$-Federer on $U$\/}  if
for all $\x\in \supp\,\mu\cap U$ one has
\begin{equation*}
\label{eq: Federer}
\frac{\mu\big(B(\x,3r)\big)}{\mu\big(B(\x,r)\big)} < D\,
\end{equation*}
whenever $B(\x,3r)\subset U$.

Say that $\mu$ is
         {\em nonplanar\/} if $\mu(\mathcal{L})=0$ for any
affine hyperplane $\mathcal{L}$ of $\R^n$.
For an affine
subspace
$\mathcal{L} \subset \R^n$
we denote by
$d_{\mathcal{L}}(\x)$ the (Euclidean) distance from $\x$ to
$\mathcal{L}$, and let
$$
\mathcal{L}^{(\varepsilon)} \stackrel{\mathrm{def}}{=} \{\x \in\R^n :
d_{\mathcal{L}}(\x) <
\varepsilon\}\,.
$$
If  $B \subset \R^n$ with $\mu(B)>0$ and $f$ is a real-valued function on
$\R^n$, let
$$\Vert f \Vert_{\mu, B} \stackrel{\mathrm{def}}{=}\sup_{x \in B \cap
\supp\, \mu} |f(x)|\,.
$$

Given $C$, $\alpha > 0$ and  an open subset $U$ of $\R^n$,   say that
$\mu$ is {\em $(C,\alpha)$-decaying on $U$} if for
any non-empty open ball $B \subset U$ centered in $\supp\,\mu$,
any affine hyperplane $\mathcal{L} \subset
\R^n$, and any $\varepsilon >0$
one has
\begin{equation*}
\label{eq: defn decaying}
{\mu \left( B \cap \mathcal{L}^{(\varepsilon)} \right) }
\le C
\left(
\frac{\varepsilon}{\Vert d_{\mathcal{L}}\Vert_{\mu, B}
}
\right)^{\alpha}{\mu(B)}\,.
\end{equation*}
Finally, let us say that
$\mu$ is {\em friendly\/} if it is nonplanar, and for $\mu$-a.e.\
$x\in \R^n$  there
exist a neighborhood
$U$ of
$x$ and positive $C, \alpha, D > 0$ such that $\mu$ is
$D$-Federer and $(C, \alpha)$ decaying on $U$.

The class of friendly measures is rather large; its properties are
discussed in \cite{friendly} and examples are given in \cite{friendly, bad,
Urbanski, Urbanski-Str}. For example, it is essentially proved in \cite{KM}
(see \cite[Lemma 7.1 and Propositions 7.2, 7.3]{friendly}) that
the natural measure on  a nondegenerate manifold
obtained by pushing forward Lebesgue measure on $\R^n$  is friendly.
Thus our  main result (Theorem \ref{cor: friendly})  supersedes
entries (a--e) in the above list. Further, there are many more possibilities
of interesting choices of sets which can support friendly measures.
A particularly nice choice is given by
limit sets of finite irreducible systems of contracting similarities
(or, more generally, self-conformal maps of $\R^n$) satisfying the
open set condition,
see \cite[\S8]{friendly} and \cite{Urbanski} for more detail.

We can now state the main measure estimate from which Theorem
\ref{cor: friendly}
will easily follow.

\begin{thm}
\label{thm: main} Let  $\mu$ be a friendly measure on $\R^n$. Then for
$\mu$-a.e.\ $\z\in\R^n$ there is a ball  $B$ centered at $\z$ and
positive $\tilde C,\alpha$ with the following property: for any $\vr$
as in \equ{eq: setting r}
there exists $t_0 > 0
$ such that for all $t > t_0$ and all $\varepsilon > 0$ one has
\begin{equation}
\label{eq: quant nondiv}
\mu\big(\{\x\in B : g^{(\vr)}_{t}\bar{\tau}(\x)
\notin K_{{\varepsilon}}\}\big)  < \tilde C\varepsilon^\alpha \mu(B)\,.
\end{equation}
\end{thm}

\begin{proof}[Proof of  Theorem \ref{cor: friendly}
assuming Theorem \ref{thm: main}] Take $B$, $\tilde C$, $\alpha$
as in Theorem \ref{thm: main},  given $\delta > 0$ choose $\varepsilon > 0$
with $\tilde C\varepsilon^\alpha \mu(B) <\delta$, and let $\{t_k\}$
be any unbounded
subsequence of $\log\mathcal{T}\cap (t_0,\infty)$.
    Then \equ{eq: ball} becomes an immediate consequence of \equ{eq:
quant nondiv}, and an application of Proposition  \ref{prop: measure
criterion}
finishes the proof.
\end{proof}


\ignore{

Thus the following `dynamical restatement' is an equivalent version
of Theorem \ref{cor: friendly}:

\begin{thm}\label{cor: friendly dynamical}
If $\mu$ is a friendly measure on $\R^n$ and $\{t_k\}\subset \R_+$ is
unbounded, then
the trajectory
$\{g^{(\vr)}_{t_k} \bar{\tau}(\x):k\in\N\}\subset \ggm
$
    is not divergent for $\mu$-a.e.\ $\x\in\R^n$.
\end{thm}
}

\section{A quantitative nondivergence estimate}

\label{quant}

We will derive Theorem \ref{thm: main}  from a quantitative nondivergence
result. To state it we need some additional definitions.

\medskip

Let $f : \R^n \to \R$.   Given $C$, $\alpha > 0$, $U \subset \R^n$
and a measure $\mu$ on
$\R^n$, say that $f$ is
         {\em $(C,\alpha)$-good on  $U$ with respect to
$\mu$\/}
         if for any ball $B \subset U$ centered in $\supp\,\mu$
and any
$\vre > 0$ one has
\begin{equation*}
\label{def-good}
{\mu\big(\{y\in B : |f(y)| < \vre\}\big)} \le C
\left(\frac{\varepsilon}{\Vert f\Vert_{\mu, B}}\right)^\alpha{\mu(B)}\,.
\end{equation*}

We refer the reader to \cite{KM, friendly} for various properties and examples.
We are going to need two elementary observations, which we state
below for convenience.

\begin{lem}
\label{lem: good} Let $U\subset \R^n$ be open, $C,\alpha > 0$, $\mu$
a measure on
$\R^n$.

\begin{itemize}
\item[(a)] \cite[Lemma 4.2]{friendly} $\mu$ is  $(C,\alpha)$-decaying
on $U$ if and only if any affine
function (equivalently, any function of the form $\,d_{\mathcal{L}}$, where
$\mathcal{L}$ is an affine hyperplane) is
$(C,\alpha)$-good on  $U$.
\item[(b)] \cite[Lemma 4.1]{friendly} If
$f_1,\dots,f_k$
are $(C,\alpha)$-good on
$U$ w.\,r.\,t.\
$\mu$, then the function  $\x\mapsto \|\vf(\x)\|$,
where $\vf = (f_1, \dots,f_k)$
and $\|\cdot\|$ is the Euclidean norm on $\R^k$, is
$(k^{\alpha/2}C,\alpha)$-good
on
$U$ w.\,r.\,t.\ $
\mu$.
\end{itemize}
\end{lem}

\ignore{
Comparing \equ{eq: defn decaying} and \equ{def-good}, one can immediately
observe, as was done in \cite[Lemma 4.2]{friendly}, that
$\mu$ is  $(C,\alpha)$-decaying on $U$ if and only if any affine
function (equivalently, any function of the form $d_{\mathcal{L}}$ where
$\mathcal{L}$ is an affine hyperplane) is
$(C,\alpha)$-good on  $U$.

In particular, we will make use of the following elementary observation,
stated in  \cite{friendly} as Lemma 4.1:
\begin{equation}
\begin{aligned}
      \label{eq: lemma}
\text{if
$f_1,\dots,f_k$
are $(C,\alpha)$-good on
$U$ with respect to
$\mu$, then} \\
\text{the function $\x\mapsto \|\vf(\x)\|$ is $(k^{\alpha/2}C,\alpha)$-good
on
$U$ w.\,r.\,t.\,}
\mu&,
\end{aligned}
\end{equation}
where $\vf = (f_1, \dots,f_k)$
and $\|\cdot\|$ is the Euclidean norm on $\R^k$.

\medskip
}

Let
$$
\mathcal{W}\df\text{ the set of nonzero rational
subspaces of } \R^{n+1}\,.
$$
For $V
\in \mathcal{W}$ and $g \in G$, let
$$
\ell_V(g) \stackrel{\mathrm{def}}{=} \Vert g (\vv_1 \wedge \cdots \wedge
\vv_k ) \Vert\,,
$$
where $\{\vv_1, \ldots, \vv_k\}$ is a generating set for $\Z^{n+1} \cap
V$ and $\Vert \cdot \Vert$ is the extension of the Euclidean norm from
$\R^{n+1}$ to its exterior algebra; note
that $
\ell_V(g)$ does not depend on the choice of  $\{\vv_i\}$.


\medskip

We will use the following estimate, which is a special case of
\cite[Theorem 4.3]{friendly}:

\begin{thm}
\label{thm: friendly nondivergence}
Given $n\in\N$ and 
positive constants $C,D,\alpha$,
         there exists $\tilde C = \tilde C(n,C,D,\alpha) > 0$
with the following property.
Suppose   $\mu$  is  $D$-Federer on an open subset $U$ of $\R^n$,
$h$ is a continuous map $U
\to G$, $0 < \rho \le 1$,
$\z\in U\cap\,\supp\,\mu$,
and   $B
   = B(\z,r)$
is a ball 
such that
denote $B(\x,cr)$
$B(\z,3^nr)\subset U$, and
that for each $V \in \mathcal{W}$,
\begin{itemize}
\item[(i)]
the function
$ \ell_V\compose {h}$  is $(C,\alpha)$-good on $B(\z,3^nr)$
with respect to
$\mu$,
           \end{itemize}
and
\begin{itemize} \item[(ii)]
\label{item: attain rho}
$\Vert \ell_V \compose {h} \Vert_{\mu,B} \geq \rho$.
          \end{itemize}
Then for any $0<
\varepsilon
\leq \rho$,
\begin{equation}
\label{eq: friendly nondiv}
{\mu\big(\big\{\x \in B: \pi\big({h}(\x)\big) \notin K_{\varepsilon}
\big\}\big)}\le \tilde C
(\varepsilon/\rho)^{\alpha}{\mu(B)} \,.
\end{equation}
\end{thm}

\begin{proof}[Proof of Theorem \ref{thm: main}] Recall that we are
given a friendly measure $\mu$.
For $\mu$-almost every $\z \in \R^n$, choose a neighborhood $U$ of
$\z$, positive constants $C', D, \alpha$ such that $\mu$ is
$D$-Federer and $(C', \alpha)$-decaying on $U$, and a ball $B = B(\z,r)$
centered at $\z$ such that $B(\z,3^nr)$ is contained in $U$.
Clearly the desired estimate \equ{eq: quant nondiv} will coincide with
   \equ{eq: friendly nondiv} if one  takes $\rho = 1$ and
lets $h = h_{\vr, t}$ where the latter is
defined by
$$h_{\vr, t}(\x) \df  g^{(\vr)}_t \tau(\x)\,.$$
Therefore it suffices to verify the assumptions of
Theorem \ref{thm: friendly nondivergence} for the above choice of
$\rho$ and $h$.
This is done below in Lemmas \ref{lem: i} and \ref{lem: ii}.
\end{proof}

\begin{lem}
\label{lem: i} Suppose that  $\mu$   is
   $(C', \alpha)$-decaying on an open $U\subset \R^n$. Then
for
any $\,V \in \mathcal{W}$,
any $\,\vr$ as in \equ{eq: setting  r} and any $\,t
\geq 0$,
\begin{equation*}
\label{eq: i}
\text{the function
$ \ell_V\compose {h_{\vr, t}}$  is $(C,\alpha)$-good
on $U$
with respect to }
\mu \,,
\end{equation*}
where $C = (n+1)^{\alpha / 2 } C'$.
\end{lem}

\begin{lem}
\label{lem: ii} Suppose that $\mu$  is nonplanar. Then for
any  $\,\vr$ as in \equ{eq: setting  r} and any ball $B$ with
$\mu(B)>0$ there is
$\,t_0= t_0(\mu,\vr,
B)$ such that for
any $\,V \in \mathcal{W}$ and any $\,t \geq t_0$ one has
\begin{equation*}
\label{eq: ii}
\Vert \ell_V \compose {h_{\vr, t}} \Vert_{\mu,B} \geq 1 \,.
\end{equation*}
\end{lem}


The proof of both lemmas hinges on
a computation of the $h_{\vr, t}(\x)$-action on the exterior powers
of $\R^{n+1}$, as in \cite{KM} or
\cite{friendly}. We include the argument for the sake of completeness.
For the remainder of the section, to simplify notation we will write
$g_t$ instead of $g^{(\vr)}_t$ and $h_{t}$  instead of $h_{\vr, t}$.

\medskip

Denote by $V_0$ the subspace
\begin{equation*}
V_0 \df \left\{ (x_1,\dots,x_{n+1}) : x_{n+1} = 0 \right\}
\end{equation*}
of $\R^{n+1}$,  and let $\ve_0 \df (0,\dots,0,1)$ be a vector orthonormal
to $V_0$. Note that for any $\x\in\R^n$,
\begin{equation}
\label{eq: action}
\tau(\x)\text{ acts trivially on $V_0$ \quad
and \quad
}\tau(\x)\ve_0 = \ve_0 + \x\,,
\end{equation}
where we with some abuse of notation identified  $V_0$ with $\R^{n}$.

Now suppose that $V$ is a  $k$-dimensional subspace of $\R^{ n + 1
}$, $k\ge 1$,
spanned by  integer vectors
$\mathbf v_1, \dots , \mathbf v_k
$, and denote $\mathbf v_1 \wedge \dots  \wedge  \mathbf v_k
$ by $\vw$.
By applying Gaussian elimination over the integers to
$\{\mathbf v_1, \dots , \mathbf v_k
\}$ one can
write $\vw$ in the form
\begin{equation}
\label{eq: elimination}
\vw = \vw_0\wedge (q \ve_0 - \p)\,,
\end{equation}
where $q\in\Z$, $\p\in V_0(\Z)$ and $\vw_0\in
\bigwedge^{k-1}\big(V_0(\Z)\big)$.
Using \equ{eq: action} and \equ{eq: elimination}, one writes
\begin{equation*}
\label{eq: tau action on w}
\tau(\x)\vw = \vw_0 \wedge  \big(q\ve_0 + q\x - \p \big) =
\vw_0 \wedge  (q \x - \p)  + q\vw_0 \wedge \ve_0 \,,
\end{equation*}
and hence
\begin{equation}
\label{eq: g action on w}
h_t(\x)\vw =
g_t \big(\vw_0 \wedge  (q \x - \p)\big) + q g_t(\vw_0 \wedge \ve_0) \,.
\end{equation}
Note that the two summands
in (\ref{eq: g action on w}) are orthogonal,  therefore
\begin{equation}
\begin{aligned}
      \label{eq: to explain}
\big(\ell_V \circ h_t(\x)\big)^2 &= \|h_t(\x)\vw\|^2  \\&=
q^2\| g_t(\vw_0 \wedge \ve_0)\|^2 +
\|g_t \big(\vw_0 \wedge  (q \x - \p)\big)\|^2
\,.
\end{aligned}
\end{equation}

\medskip

Now we can return to the lemmas.

\begin{proof}[Proof of Lemma \ref{lem: i}] Write the second summand
in (\ref{eq: to explain})
in the form
$$\|g_t \vw_0 \wedge  g_t(q \x - \p)\|^2 = \|g_t\vw_0\|^2
d_{g_t\mathcal{P}} \big(  g_t
( q\x - \p) \big)^2\,,
$$
where  $\mathcal{P}$  stands for the linear subspace of $V_0$ corresponding
to $\vw_0$. For any $t$, $q$ and $\p$, the function $\x\mapsto
d_{g_t\mathcal{P} }^2\big(  g_t
( q\x - \p) \big)$ is the sum of squares of at
most $n$ affine functions, each of which is $(C',\alpha)$-good on
$U$ with respect to
$\mu$ in view
of Lemma \ref{lem: good}(a)  and the assumption of Lemma  \ref{lem:
i}. Therefore,
by  Lemma \ref{lem: good}(b), the  function $\ell_V \circ h_t$
is $\big((n+1)^{\alpha / 2 } C',\alpha\big)$-good
on $U$
with respect to
$\mu$.
\end{proof}

\begin{proof}[Proof of Lemma \ref{lem: ii}]
Let us denote by $e^{\gamma t}$ the smallest eigenvalue of the induced action
of $g_t$ on $\bigwedge^{k}(V_0)$ (here $\gamma > 0$ depends on $\vr$).
If $q = 0$, in view of \equ{eq: elimination} we have
$\ww\in \bigwedge^{k}\big(V_0(\Z)\big)$, hence
$$
   \ell_V \circ h_t(\x) = \|g_t \tau(\x)\vw\| = \|g_t \vw\| \ge
e^{\gamma t} \| \vw\| \ge 1
$$
for all $t \ge 0$. Thus the conclusion of the lemma holds  with e.g.\ $t_0=0$.
Otherwise, using (\ref{eq: g action on w}),
one can write
\begin{equation*}
\begin{aligned}
   \ell_V \circ h_t(\x) &\ge
\|g_t \big(\vw_0 \wedge  (q \x - \p)\big)\|
\ge e^{\gamma t}\|\vw_0 \wedge  (q \x - \p)\|\\
& = e^{\gamma t}\|\vw_0\|  d_{\mathcal{P}}
( q\x - \p) = |q|e^{\gamma t}\|\vw_0\|  d_{(\mathcal{P} + \p/q)}
( \x) \ge e^{\gamma t}  d_{(\mathcal{P} + \p/q)}( \x)
\end{aligned}
\end{equation*}
(the last inequality holds since both $q$ and all the coordinates
of $\vw_0$ are integers).

   If $B$ is a ball with
$\mu(B) >0$, an  easy compactness argument using the assumption
that
$\mu$ is nonplanar shows the existence of
$c = c(B)>0$ such that
$\Vert d_{\mathcal{L}}\Vert_{\mu, B} \geq c\,$ for
         any proper affine  subspace $\mathcal{L}$ of $\R^n$.
Hence
\begin{equation*}
   \Vert\ell_V \circ h_t\Vert_{\mu, B} \ge c|q|e^{\gamma t}
\|\vw_0\|\ge ce^{\gamma t}
\,,
\end{equation*}
and the conclusion of the lemma holds with $t_0 =
\frac1\gamma \log\frac1c$.
\end{proof}

\section{Constructing singular vectors on submanifolds}

\label{cassels}

In this section we adapt the methods of the paper \cite{cassels}
and
exhibit points with divergent trajectories on certain
proper subsets of $G/\Gamma$.


\begin{thm}
\label{thm: like cassels}
Let $\{g_t : t \in \R\}$ be a one-parameter subgroup of $G$.
Suppose $X$ is a closed
subset of $G$, and $\{X_i: i
\in \N\}$
and $\{X'_j: j
\in \N\}$ are two lists of 
subsets of $X$,
such that for some
strictly decreasing continuous $\psi: \R_+ \to \R_+$, the
following conditions are satisfied.

\begin{enumerate}
\item{\bf Density.} For every $j$,
$X_j = \cl{X_j \cap \bigcup_{i \neq j} X_i}.$

\item{\bf Transversality I.} For every $i \neq j$,
$X_i  = \cl{X_i \sm X_j}.$

\item{\bf Transversality II.} For any $i, j$,
$X_i = \cl{X_i \sm X'_j}.$

\item{\bf Local Uniformity w.r.t.\ $\{K_{\psi(t)}\}$.} For every
$i$ and every $\,x \in X_i$
there is a neighborhood ${U}$ in $G$ and $t_0$ such that for
all $t \geq t_0$ and all $z \in {U} \cap X_i$,
$$g_t \pi(z) \notin K_{\psi(t)}.
$$
\end{enumerate}

Then there is  $x \in X \sm
\left( \bigcup_i X_i \cup \bigcup_j X'_j \right)
$ and $t_0>0$ such that
for all $t \geq t_0$, $g_t \pi(x) \notin K_{\psi(t)}.$
\end{thm}

\begin{proof}
Equip $X$ with the relative topology inherited from $G$, and write
$$K(t)=K_{\psi(t)}.
$$

We will construct inductively
a sequence of open sets with compact closure $\Omega_0, \Omega_1, \Omega_2,
\ldots$ in $X$, an increasing sequence of
indices $i_1, i_2, \ldots$, and an increasing sequence of times $T_0,
T_1, \ldots, $ such that the following hold for $k=1,2, \ldots$:
\begin{itemize}
\item[a.]
$\cl{\Omega_{k}} \subset \Omega_{k-1}$.
\item[b.]
For every $j < i_k$, $X_j \cap \Omega_k = \varnothing$ and $X'_j \cap
\Omega_k = \varnothing$.
\item[c.]
$X_{i_k} \cap \Omega_k$ is nonempty and
for every $z \in X_{i_k} \cap \Omega_k$ and every $t \geq T_k$ we
have $g_t\pi(z) \notin K(t).$
\end{itemize}
We will also have for $k=2,3, \ldots$:
\begin{itemize}
\item[d.]
For every $z \in \Omega_k$ and every $t \in [T_{k-1}, T_k]$,
$g_t\pi(z) \notin K(t)$.

\end{itemize}

To see that such sequences suffice, note that by condition a,
$\bigcap_k \Omega_k$ is nonempty, and for $z \in \bigcap \Omega_k$ we
have by condition b that $z \notin \bigcup_i X_i \cup \bigcup_j X'_j$
and by condition d that $g_t \pi(z) \notin K(t)$ for $t \geq T_1$.

\medskip

Now let us construct the sequences inductively.
Choose $T_0=0$, $i_1=1$. Let $x
\in X_1$ and using the local uniformity hypothesis, let $\Omega_1$ be
a small enough open neighborhood of
$x$, and $T_1$ large enough, so that for all $z \in X_1 \cap \Omega_1$
and all $t \geq T_1,$ we have $g_t \pi(z) \notin K(t).$ Now letting
$T_0$ be arbitrary and $\Omega_0$ be any open set containing
$\cl{\Omega_1}$, we see that a, b and c hold for $k=1$.

Suppose we have chosen $i_s, \Omega_s, T_s$ for $s=1, \ldots, k$. By
the density condition there are $\ell \neq i_k$ such that
$$X_{\ell}
\cap \Omega_k \cap X_{i_k} \neq \varnothing.$$
Choose for $i_{k+1}$ any such $\ell$. Note that $i_{k+1} > i_k$ by
b. Let $x \in X_{i_k} \cap \Omega_k \cap X_{i_{k+1}}.$ By the local
uniformity assumption, there is a small enough open neighborhood
${U}$ of $x$ and a large enough $T_{k+1}$ such that for all $z
\in {U} \cap X_{i_{k+1}}$ and all $t \geq T_{k+1}$, $g_t
\pi(z) \notin K(t).$ In addition let ${U}$ be small enough so
that $\cl{{U}} \subset \Omega_k$. Since $x \in X_{i_k}$, $g_t
\pi(x)  \notin K(t)$ for $t \in [T_k, T_{k+1}]$, hence by continuity
of $\psi$ and the action, there is a small enough neighborhood
$\til{\Omega}$ of $x$
contained in ${U}$ so that
$$ z \in \til{\Omega}, \ t \in [T_{k-1}, T_k] \ \ \ \Longrightarrow \
\ \ g_t\pi(z) \notin K(t).$$
Now we can define $\Omega_{k+1}$ by
$$\Omega_{k+1}=\til{\Omega} \sm  \bigcup_{j < i_{k+1}}
\left(X_j \cup X'_j\right).
$$

We now verify that $i_{k+1}, \Omega_{k+1}, T_{k+1}$
satisfy the required conditions. Condition a holds by our choice of
${U}.$ Condition b follows from the definition of
$\Omega_{k+1}$.  In
condition c, $\Omega_{k+1} \cap X_{i_{k+1}} \neq \varnothing$ because $x
\in \til{\Omega} \cap X_{i_{k+1}}$, and because of the
transversality assumptions. The second assertion in condition c holds
because of the choice of $T_{k+1}$ and ${U}$. Condition d holds because
of the choice of $\til{\Omega}$.
\end{proof}

We now derive a consequence of this theorem. This requires some
notation.
\ignore{
Given $\vr=(r_1, \dots, r_n)$ as in \equ{eq: setting r}, let $\rho_2 =
\rho_2(\vr)$
be the second largest of $r_1, \ldots, r_n$, that is,
$$\rho_2 = r_{i_2}, \ \ \ \mathrm{where \ } r_{i_1} \geq r_{i_2} \geq
\cdots \geq r_{i_n}, \ \ \{1, \ldots, n\} = \{i_1, \ldots, i_n\}.$$
}
\ignore{
Let $U \subset \R^d$ be open and $\ff: U \to
\R^n$ be a $C^{\ell}$ map for some $\ell \in \N$. Then $\ff$ is called
{\em non-degenerate} if for almost every $u \in
U$, the partial derivatives of
$\ff$ up to order $\ell$ at $u$ span $\R^n.$ \combarak{This is
equivalent to the usual definition, isn't it?} Similarly, ${M} \subset
\R^n$ is called {\em non-degenerate} if there is a non-degenerate
$\ff: U \to \R^n$ as above with ${M}=\ff(U).$
}
Fix $\vr$, let ${{M}}$ be a submanifold of  $\R^n$, and choose $1 \leq
k < \ell \leq n$. For $\vv \in \Z^n$ and $s \in \Q$ let
$$L_{\vv}(s)=\{\x \in {M}: \langle \x, \vv \rangle=s\}.$$
Let $\E_k, \E_{\ell}$ be the $k$-th and $\ell$-th standard basis vectors
and let $\{X_i\}$ be a list of the distinct connected
components of the sets $\{L_{\E_{r}}(s): r \in \{k, \ell\}, \, s \in
\Q\}$. Also
let $\{X'_j\}$ be a list of the distinct connected components of the
sets $\{L_{\vv}(s): \vv \in \Z^n, \, s \in \Q\}$ which are not in the
list $\{X_i\}$. In case
$\{X'_j\} = \{X'_1, \ldots, X'_m\}$ happens to be a finite list, we
put $X'_j = \varnothing$ for $j>m$, i.e., we may assume that
$\{X'_j\}$ is indexed by $\N$ as well.

\begin{cor}
\label{cor: main, dim2}
Suppose $\vr, \, {M}, \, k, \, \ell, \, \{X_i\},\, \{X'_j\}$ are as above, and
that the density and two
transversality assumptions of Theorem \ref{thm: like cassels} are
satisfied. Suppose also that
$\delta: \R_+ \to \R_+$ is a function such that
   $$T^{\rho} \delta(T) \to_{T \to \infty} \infty, \ \ \ \mathrm{where \ }
\rho =
\min \{r_{k}, r_{\ell}\}/n
.$$

Then there is a totally irrational $\x \in
{M}$ and $T_0$ such that for all $T \geq T_0$ there is $\pp
\in \Z^n, \, q \in \N$ satisfying
\equ{eq: r-T-singular}  for $\delta = \delta(T).$

\end{cor}

\begin{proof}
Set
$\psi(t) \df \delta(e^t)$, so that
\begin{equation}
\label{eq: speed}
e^{\rho t} \psi(t) \to \infty.
\end{equation}
Suppose without loss of generality that $\psi$ is a decreasing function.
Using Lemma \ref{lem: correspondence},
it suffices to show that there are $x = \tau(\x) \in
\tau({M})$ and $t_0$ such that $\x$ is totally irrational and for all $t \geq
t_0$,
$$g^{(\vr)}_t \pi(x) \notin K_{\psi(t)}.
$$

To this end, we apply Theorem \ref{thm: like cassels} to the lists
$\{X_i\}, \, \{X'_j\}$ (which we identify with their images under
$\tau$, and thus consider them as subsets of $G$). Since we have
assumed the density and
two transversality conditions, we need only check the locally uniform
escape condition.

To
verify the condition of local uniformity w.r.t.\ $\{K_{\psi(t)}\}$, fix
$X_i$ so that (possibly after exchanging $k$ and $\ell$), for some $p/q
\in \Q$ and for all $\z=(z_1, \ldots, z_n) \in
X_i$ we have $z_{k} = p/q$.
Let
$$W  = \{(w_1, \ldots, w_{n+1}) \in \R^{n+1} : qw_{k}+pw_{n+1} = 0\}.$$
This is a rational linear subspace of $\R^{n+1}$ of dimension
$n$, with the property that for all $\z \in X_i$, all vectors in
$\tau(\z) \cdot W$ have their
$k$-th coordinate equal to zero. It follows that for any bounded
neighborhood ${U}$ intersecting $X_i$ there is a
constant $C$, such that for all $\z \in {U} \cap X_i$ we have
$$\ell_W\left(g^{(\vr)}_t \tau(\z)\right) \leq C e^{-r_{k} t}.$$
Using Minkowski's convex body theorem, we find a constant $C'$ so that
for all $t$ and all $\z \in {U} \cap X_i$,
$g^{(\vr)}_t \tau(\z) \, \left(W \cap \Z^{n+1}\right)$ contains a
non-zero vector of
length at most
$$C' e^{-r_{k}t/n} \leq C' e^{-\rho t}.$$
Now from \equ{eq: speed} if follows that there is $t_0$ such that for
all $t \geq t_0$, the length of such a vector is less
than $\psi(t)$. This concludes the proof.
\end{proof}

\ignore{\begin{cor}
\label{cor: analytic manifolds}
Fix $\vr$, and let ${M} \subset \R^n$ be an analytic
submanifold, of
dimension at least $2$, which is not contained in any proper rational affine
hyperplane.
Then ${M} \cap {\Sing}(\vr)$ contains a totally irrational vector.
\end{cor}}

\begin{proof}[Proof of Theorem \ref{thm: main, dim2}]
For $\z \in {M}$ let $T_\z{M} \subset \R^n$ be the tangent space to
${M}$ at $\z$. Since $\dim ({M}) \geq 2$, there are distinct indices $k,
\ell \in \{1,
\ldots, n\}$ and an open $V \subset {M}$ such that
  if $P$ is the projection
$$P(\x) = (x_{k}, x_{\ell}) \ \ \ \mathrm{for \ } \x = (x_1, \ldots,
x_n),$$
then the
derivative $D_{\z} \left(P|_{{M}}\right): T_\z M \to \R^2$ is
surjective for all $\z \in V$.
With no loss of generality, replace ${M}$ with $V$  and define
$\{X_i\}, \, \{X'_j\}$ as in the paragraph preceding Corollary
\ref{cor: main, dim2}. In view of the Corollary, it remain to check
that the density and two transversality hypotheses hold for $\{X_i\},
\, \{X'_j\}$.

Since these hypotheses hold for horizontal and vertical
lines in an open subset of the plane, and since the $X_i$'s are the
pre-images of these lines under $P$, we see that the density and the first
transversality hypothesis hold. Now for $i$ and $j$, suppose that $X_i
\cap X'_j \neq \varnothing$ (otherwise there is nothing to
prove). Then $X_i$ and $X_i \cap X'_j$ are connected analytic submanifolds of
${M}$. Suppose if possible that $X_i \cap X'_j$ contains an open
subset of $X_i$. Then, since they are analytic and $X_i$ is connected, $X_i
\subset X'_j$. Since ${M}$ is not contained in a rational affine
hyperplane, both $X_i$ and $X'_j$ are submanifolds of ${M}$ of
codimension one, and since $X'_j$ is also connected we must have $X_i
= X'_j$, contrary to the
construction. This implies
that $X_i = \cl{X_i \sm X'_j}$, as required.
\end{proof}

\begin{remark} \rm
The hypotheses of Theorem \ref{thm: main, dim2} are clearly satisfied
when ${M}$ is an analytic
nondegenerate submanifold of dimension at least $2$, and also when
${M}$ is an affine subspace of dimension at least $2$ not contained in
any rational affine hyperplane.
\end{remark}

\begin{remark} \rm
The proof of Theorem 1.2 actually shows that $M \cap \Sing
(\vr)$ contains {\em uncountably many} totally irrational
vectors. Indeed, given any countable subset  $A= \{\z_1, \z_2, \ldots
\}$ of $M$, replace the list $\{X'_j\}$ with the list $\big
\{X'_1, \{\z_1\}, X'_2, \{\z_2\}, \ldots \big\}$. Applying the same
argument yields an element of $M \cap \Sing(\vr)$ which is totally
irrational and does not belong to $A$.
\end{remark}



\end{document}